\documentclass[envcountsame]{mmnp}
\usepackage[tbtags]{amsmath}
\usepackage{amssymb}
\usepackage{graphicx}
\usepackage[latin1]{inputenc}

\idline{Vol. 7, No. 2}{32}
\doi{10.1051/mmnp/20127102}

\numberwithin{equation}{section}

\begin{document}

\title{Heat transfer in a medium in which many small particles are embedded}


\author{A. G. Ramm\inst{1}   \thanks{\email {ramm@math.ksu.edu}} }

\vspace{0.5cm}

\institute{\inst{1} Department of Mathematics \\ Kansas State
University, Manhattan, KS 66506-2602, USA }


\abstract{The heat equation is considered in the complex system
consisting of many small bodies (particles) embedded in a given
material. On the surfaces of the small bodies a Newton-type boundary
condition is imposed. An equation for the limiting field is derived
when the characteristic size $a$ of the small bodies tends to zero,
their total number $\mathcal{N}(a)$ tends to infinity at a suitable
rate, and the distance $d = d(a)$ between neighboring small bodies
tends to zero $a << d$. No periodicity is assumed about the
distribution of the small bodies. }





\titlerunning{Heat transfer in a medium }

\maketitle


\section{Introduction}\label{Introduction} Let many small
bodies (particles) $\mathcal{D}_m$, $1 \leq m \leq M$, be
distributed in a bounded domain $\mathcal{D} \subset \mathbb{R}^3$,
diam$\mathcal{D}_m = 2a$. The small bodies are distributed according
to the law
\begin{equation}\label{eq:1} \mathcal{N}(\Delta) = \frac{1}{a^{2 -
\kappa}}\int_{\Delta} N(x)dx[1 + o(1)] , a \rightarrow 0. \end{equation}
Here $\Delta \subset \mathcal{D}$ is an arbitrary open subdomain of
$\mathcal{D}$, $\kappa \in (0, 1)$ is a constant, $N(x) \geq 0$ is a
continuous function, and $\mathcal{N}(\Delta)$ is the number of the
small bodies $\mathcal{D}_m$ in $\Delta$. Let us assume that the
boundaries $S_m$ of the bodies $\mathcal{D}_m$ are $C^{2}-$smooth.

The heat equation can be stated as follows:
\begin{equation}\label{eq:2} u_t = \nabla^2u + f(x) \, \textrm{ in }
\, \mathbb{R}^3 \setminus \displaystyle\bigcup^{M}_{m =
1}\mathcal{D}_m := \Omega, \, u |_{t = 0} = 0, \end{equation}
\begin{equation}\label{eq:3} u_N = \zeta_m u \textrm{ on } \mathcal{S}_m,
1 \leq m \leq M. \end{equation} Here $N$ is the outer unit normal to
$\mathcal{S} := \displaystyle\bigcup^{M}_{m = 1} \mathcal{S}_m$,
$\zeta_m = \displaystyle\frac{h(x_m)}{a^{\kappa}}, x_m \in
\mathcal{D}_m, 1 \leq m
\leq M$, where $h(x)$ is a continuous function in $\mathcal{D}$.\\

The aim of this paper is to establish a method for solving both
theoretically and numerically the many-body heat transfer problem,
and to derive an equation for the effective (self-consistent) field
in the medium in which many small bodies (particles) are embedded.

There is a large literature on homogenization theory which deals
with similar problems, see \cite{JKO}-\cite{BLP} and the literature
therein. In most cases it is assumed in the homogenization
literature that the coefficients of the equations are periodic, and
there is a small parameter in these coefficients, that the related
problems in a "periodic cell" have discrete spectrum and are
selfadjoint, that the elliptic estimates, such as Rellich and
Poincare inequalities hold, etc. Our method, which was used in
\cite{R476}-\cite{R610}, differs in many respects from the
homogenization theory methods developed in the literature:
periodicity assumption is not used, the spectrum of the related
elliptic problems that we use is continuous, the problems may be
non-selfadjoint, and our justification of the limiting equation for
the effective field is based on a new technique, namely, on
convergence results for collocation methods (see, \cite{R595}). The
published homogenization techniques are not directly applicable to
our problem because of the lack of periodicity.

The main results of this paper are:

1) Derivation of the linear algebraic system (25) for finding the
Laplace transform of the solution to problem (2)-(3),

2) Derivation of the equation (29) for the limiting effective field
in the medium when the size of the small bodies tends to zero while
the total number of these bodies tends to infinity according to the
distribution law (1),

3) Derivation of the formula (31) for the average temperature in the
limiting medium.

In Section 2 a derivation of these results is given, and in Section
3 proofs of some lemmas are given.

\section{Derivation of the equation for the limiting  effective field}
\label{Derivation of the equation for the limiting  effective field}
Denote $$\mathcal{U} := \mathcal{U}(x, \lambda) =
\displaystyle\int^{\infty}_{0}e^{-\lambda t}u(x, t)dt.$$ Then,
(\ref{eq:2}) - (\ref{eq:3}) imply \begin{equation}\label{eq:4}
-\nabla^2\mathcal{U} + \lambda\mathcal{U} = \lambda^{-1}f(x)
\,\,\textrm{ in } \,\, \Omega, \end{equation}
\begin{equation}\label{eq:5} \mathcal{U}_N = \zeta_m\mathcal{U} \,\, \textrm{
on }\,\, \mathcal{S}_m, \,\, 1 \leq m \leq M. \end{equation}
Let
\begin{equation}\label{eq:6} g(x,y) := g(x, y, \lambda) =
\frac{e^{-\sqrt{\lambda}|x - y|}}{4\pi|x - y|}, \quad
\frac{1}{\lambda}\int_{\mathbb{R}^3}g(x,y)f(y)dy = F(x, \lambda).
\end{equation} Look for the solution to (\ref{eq:4}) - (\ref{eq:5}) of the
form \begin{equation}\label{eq:7} \mathcal{U}(x, \lambda) = F(x,
\lambda) + \sum^{M}_{m = 1}\int_{\mathcal{S}_m}g(x,
s)\sigma_m(s)ds,\quad \mathcal{U}(x,
\lambda):=\mathcal{U}(x):=\mathcal{U},
\end{equation}
where $\sigma_m$ are unknown and should be found from the boundary
conditions (\ref{eq:5}). Equation (\ref{eq:4}) is satisfied by
$\mathcal{U}$ of the form (\ref{eq:7}) for any $\sigma_m$. To
satisfy (\ref{eq:5}) one has to solve equation
\begin{equation}\label{eq:8} \frac{\partial
\mathcal{U}_e(x)}{\partial N} + \frac{A_m\sigma_m - \sigma_m}{2} -
\zeta_m\mathcal{U}_e - \zeta_mT_m\sigma_m = 0 \textrm{ on }
\mathcal{S}_m, 1 \leq m \leq M. \end{equation} Here
\begin{equation}\label{eq:9} \mathcal{U}_e(x) := \mathcal{U}_{e,m}(x) :=
\mathcal{U}(x) - \int_{\mathcal{S}_m}g(x, s)\sigma_m(s)ds, \end{equation}
\begin{equation}\label{eq:10} T_m\sigma_m = \int_{\mathcal{S}_m}g(s,
s')\sigma_m(s')ds', \, A_m\sigma_m = 2\int_{\mathcal{S}_m}\frac{\partial
g(s, s')}{\partial N_S}\sigma_m(s')ds', \end{equation} and the known
formula for the outer limiting value on $\mathcal{S}_m$ of the
normal derivative of a simple layer potential was used. We now apply
the ideas and methods for solving many-body scattering problems
developed in \cite{R476} - \cite{R595}.$\\\\$ Let us call
$\mathcal{U}_{e,m}$ the effective (self-consistent) value of
$\mathcal{U}$, acting on $m$-th body. As $a \rightarrow 0$, the
dependence on $m$ disappears, since
$$\displaystyle\int_{\mathcal{S}_m}g(x, s)\sigma_m(s)ds \rightarrow 0
\,\,\textrm{ as }\,\, a \rightarrow 0.$$ One has
\begin{equation}\label{eq:11} \mathcal{U}(x, \lambda) = F(x,
\lambda) + \sum^{M}_{m = 1}g(x, x_m)Q_m + \mathcal{J}_2, \quad x_m
\in \mathcal{D}_m, \end{equation} where $$Q_m :=
\int_{\mathcal{S}_m}\sigma_m(s)ds,$$ and
\begin{equation}\label{eq:12} \mathcal{J}_2 := \sum^{M}_{m =
1}\int_{\mathcal{S}_m}[g(x, s') - g(x, x_m)]\sigma_m(s')ds', \,\,
\mathcal{J}_1 := \sum^{M}_{m = 1}g(x, x_m)Q_m.
\end{equation} We prove (in Section 2) that \begin{equation}\label{eq:13}
|\mathcal{J}_2| << |\mathcal{J}_1| \textrm{ as } a \rightarrow 0
\end{equation} provided that \begin{equation}\label{eq:14} \lim_{a
\rightarrow 0}\frac{a}{d(a)} = 0. \end{equation} If (\ref{eq:13})
holds, then problem (\ref{eq:4}) - (\ref{eq:5}) is solved
asymptotically by the formula \begin{equation}\label{eq:15}
\mathcal{U}(x, \lambda) = F(x, \lambda) + \sum^{M}_{m = 1}g(x,
x_m)Q_m, \, a \rightarrow 0,
\end{equation} provided that asymptotic formulas for $Q_m$, as $a
\rightarrow 0$, are found. To find formulas for $Q_m$, let us
integrate (\ref{eq:8}) over $\mathcal{S}_m$ and estimate the order
of the terms in the resulting equation as $a \rightarrow 0$. We get
\begin{equation}\label{eq:16} \int_{\mathcal{S}_m}\frac{\partial
\mathcal{U}_e}{\partial N}ds = \int_{\mathcal{D}_m}\nabla^2\mathcal{U}_edx
= O(a^3). \end{equation} Here we assumed that $|\nabla^2\mathcal{U}_e| =
O(1), a \rightarrow 0$. This assumption will be justified in Section
2.
\begin{equation}\label{eq:17} \int_{\mathcal{S}_m}\frac{A_m\sigma_m -
\sigma_m}{2}ds = - Q_m[1 + o(1)], \,\,\, a \rightarrow 0. \end{equation} This
relation is justified in Section 2. Furthermore,
\begin{equation}\label{eq:18} -\zeta_m\int_{\mathcal{S}_m}\mathcal{U}_eds
= -\zeta_m|\mathcal{S}_m|\mathcal{U}_e(x_m) = O(a^{2 - \kappa}), \,\,\, a
\rightarrow 0, \end{equation} where $$|\mathcal{S}_m| = O(a^2)$$ is the
surface area of $\mathcal{S}_m$. Finally,
$$-\zeta_m\int_{\mathcal{S}_m}ds\int_{\mathcal{S}_m}g(s,
s')\sigma_m(s')ds' =
-\zeta_m\int_{\mathcal{S}_m}ds'\sigma_m(s')\int_{\mathcal{S}_m}dsg(s,
s')$$ \begin{equation}\label{eq:19}
 = Q_m O(a^{1 - \kappa}), \quad a \rightarrow 0. \end{equation} Thus, the
main term of the asymptotics of $Q_m$ is
\begin{equation}\label{eq:20} Q_m =
-\zeta_m|\mathcal{S}_m|\mathcal{U}_e(x_m). \end{equation} Formulas
(\ref{eq:20}) and (\ref{eq:15}) yield \begin{equation}\label{eq:21}
\mathcal{U}(x, \lambda) = F(x, \lambda) - \sum^{M}_{m =
1}\zeta_m|\mathcal{S}_m|\mathcal{U}_e(x_m, \lambda), \end{equation}
and
\begin{equation}\label{eq:22} \mathcal{U}_e(x_m, \lambda) = F(x_m,
\lambda) - \sum^{M}_{m' \neq m, m' = 1}g(x_m,
x_{m'})\zeta_{m'}|\mathcal{S}_{m'}|\mathcal{U}_e(x_{m'}, \lambda).
\end{equation}
Denote $$\mathcal{U}_e(x_m, \lambda) := \mathcal{U}_m,$$
$$F(x_m, \lambda) := F_m, \quad g(x_m, x_{m'}) := g_{mm'},$$ and write
(\ref{eq:22}) as a linear algebraic system
\begin{equation}\label{eq:23} \mathcal{U}_m = F_m - a^{2 -
\kappa}\sum_{m' \neq m}g_{mm'}h_{m'}c_{m'}\mathcal{U}_{m'},\quad 1
\leq m \leq M,
\end{equation}
where $$h_{m'} = h(x_{m'}),\quad \zeta_{m'} =
\frac{h_{m'}}{a^{\kappa}}, \quad c_{m'} := |S_{m'}|a^{-2}. $$
Consider a partition of the bounded domain $\mathcal{D}$, in which
the small bodies are distributed, into a union of $P << M$ small
nonintersecting cubes $\Delta_p$, $1 \leq p \leq P$, of side $b >>
d$, $b = b(a) \rightarrow 0$ as $a \rightarrow 0.$ Let $x_p \in
\Delta_p$, $|\Delta_p| =$ volume of $\Delta_p$. One has $$a^{2 -
\kappa}\sum^{M}_{m' = 1, m' \neq
m}g_{mm'}h_{m'}c_{m'}\mathcal{U}_{m'} = a^{2 - \kappa}\sum^{P}_{p' =
1, p' \neq p}g_{pp'}h_{p'}c_{p'}\mathcal{U}_{p'}\sum_{x_{m'} \in
\Delta_{p'}}1 =$$
\begin{equation}\label{eq:24} = \sum_{p' \neq
p}g_{pp'}h_{p'}c_{p'}\mathcal{U}_{p'}N(x_{p'})|\Delta_{p'}|[1 + o(1)],
\quad a \rightarrow 0. \end{equation} Thus, (\ref{eq:23}) yields
\begin{equation}\label{eq:25} \mathcal{U}_p = F_p - \sum^{P}_{p' \neq p,
p' = 1}g_{pp'}h_{p'}c_{p'}N_{p'}\mathcal{U}_{p'}|\Delta_{p'}|, \quad 1
\leq p
\leq P \end{equation}
We have assumed that \begin{equation}\label{eq:26}
h_{m'} = h_{p'}[1 + o(1)],\quad c_{m'} = c_{p'}[1 + o(1)],\quad
\mathcal{U}_{m'} = \mathcal{U}_{p'}[1 + o(1)], \,\,a \rightarrow 0,
\end{equation} for $x_{m'} \in \Delta_{p'}$. This assumption is
justified if the functions $h(x), \mathcal{U}(x, \lambda)$, $$c(x) =
\displaystyle\lim_{x_{m'} \in \Delta_x, a \rightarrow
0}\displaystyle\frac{|S_{m'}|}{a^2},$$ and $N(x)$ are continuous.
The function $h(x)$ and $N(x)$ are continuous by the assumption. The
continuity of the $\mathcal{U}(x, \lambda)$ is proved in Section 2,
and the continuity of $c(x)$ is assumed. If all the small bodies are
identical, then $c(x) = c =$ const.$\\$ The sum in the right-hand
side of (\ref{eq:25}) is the Riemannian sum for the integral
\begin{equation}\label{eq:27} lim_{a \rightarrow 0}\sum^{P}_{p' = 1, p'
\neq p}g_{pp'}h_{p'}c_{p'}N(x_{p'})\mathcal{U}_{p'} =
\int_{\mathcal{D}}g(x,y)h(y)c(y)N(y)\mathcal{U}(y, \lambda)dy
\end{equation} Therefore, linear algebraic system (\ref{eq:25}) is a
collocation method for solving integral equation
\begin{equation}\label{eq:28} \mathcal{U}(x, \lambda) = F(x, \lambda) -
\int_{\mathcal{D}}g(x,y)h(y)c(y)N(y)\mathcal{U}(y, \lambda)dy.
\end{equation} Convergence of this method for solving equations with
weakly singular kernels is proved in \cite{R563}.$\\$ Applying the
operator $-\nabla^2 + \lambda$ and then taking the inverse Laplace
transform of (\ref{eq:28}) yields \begin{equation}\label{eq:29} u_t
= \Delta u + f(x) - q(x)u, \, q(x) := h(x)c(x)N(x). \end{equation}
One concludes that the limiting equation for the temperature
contains the term $q(x)u$. Thus, the embedding of many small
particles creates a distribution of source and sink terms in the
medium, the distribution of which is described by the term
$q(x)u$.$\\$ If one solves equation (\ref{eq:28}) for
$\mathcal{U}(x, \lambda)$, or linear algebraic system (\ref{eq:25})
for $\mathcal{U}_p(\lambda)$, then one can Laplace-invert
$\mathcal{U}(x, \lambda)$ for $\mathcal{U}(x, t)$. Numerical methods
for Laplace inversion from the real axis are discussed in
\cite{R198} - \cite{R569}.$\\$ If one is interested only in the
average temperature, one can use the relation
\begin{equation}\label{eq:30} \lim_{T \rightarrow
\infty}\frac{1}{T}\int^{T}_{0}u(x, t)dt = \lim_{\lambda \rightarrow
0}\lambda\mathcal{U}(x, \lambda) := \psi(x). \end{equation} Relation
(\ref{eq:30}) is proved in Section 2, which holds if the limit on
one of its sides exists. The limit on the right-hand side of
(\ref{eq:30}) can be calculated by the formula
\begin{equation}\label{eq:31} \psi(x) = (I + B)^{-1}\varphi, \,\,
\varphi = \int\frac{1}{4\pi|x - y|}f(y)dy. \end{equation} Here, $B$
is the operator $$B\psi :=
\displaystyle\int\displaystyle\frac{q(y)\psi(y)}{4\pi|x -
y|}dy,\quad q(x) := h(x)c(x)N(x).$$ From the physical point of view
the function $h(x)$ is non-positive because the flux $-\nabla u$ of
the heat flow is proportional to the temperature $u$ and is directed
along the outer normal $N$: $-u_N = h_1u$, where $h_1 = -h > 0$.
Thus, $q \leq 0$. It is proved in \cite{R203} - \cite{R227} that
zero is not an eigenvalue of the operator $-\nabla^2 + q(x)$
provided that $q(x) \geq 0$ and $q = O\big(
\displaystyle\frac{1}{|x|^{2 + \epsilon}} \big)$ as $|x| \rightarrow
\infty, \epsilon > 0.$ In our case, $q(x) = 0$ outside
$\mathcal{D}$, so the operator $(I + B)^{-1}$ exists and is bounded
in $C(\mathcal{D})$. Let us formulate the basic result we have
proved.

{\bf Theorem 1.} {\it Assume (\ref{eq:1}), (\ref{eq:14}), and $h
\leq 0$. Then, there exists the limit $\mathcal{U}(x, \lambda)$ of
$\mathcal{U}_e(x, \lambda)$ as $a \rightarrow 0$, $\mathcal{U}(x,
\lambda)$ solves equation (\ref{eq:28}), and there exists the limit
(\ref{eq:30}), where $\psi(x)$ is given by formula
(\ref{eq:31}).}
\section{Proofs of some lemmas} \label{Proofs  of some lemmas}
\begin{lmm}
Assume (\ref{eq:14}). Then relation (\ref{eq:13}) holds. \end{lmm}
\begin{proof} One has \begin{equation}\label{eq:32} \mathcal{J}_{1, m} :=
|g(x, x_m)Q| = O\bigg( \frac{|Q_m|e^{-\sqrt{x}|x - x_m|}}{4\pi|x -
x_m|} \bigg) \leq \frac{e^{-1}|Q_m}{|x - x_m|}, |x - x_m| \geq d.
\end{equation}
\begin{equation}\label{eq:33} \mathcal{J}_{2, m} :=
\int_{\mathcal{S}_m}\frac{e^{-\sqrt{\lambda}|x - x_m|}}{4\pi|x -
x_m|}\max\bigg( \sqrt{\lambda}a, \frac{a}{|x - x_m|}
\bigg)|\sigma_m(s')|ds' \leq O\bigg( \frac{|Q_m|a}{|x - x_m|^2} \bigg),
\end{equation} where $|x - x_m| \geq 2$, and the inequality
$$\displaystyle\max_{\lambda \geq
0}(\sqrt{\lambda}e^{-\sqrt{\lambda}|x - x_m|}) \leq
\displaystyle\frac{e^{-1}}{|x - x_m|}$$ was used. The $|Q_m| \neq 0$.
In fact, $\sigma_m$ keeps sign on $\mathcal{S}_m$, as follows from
equation (\ref{eq:8}) as $a \rightarrow 0$.$\\$ It follows from
(\ref{eq:32}) - (\ref{eq:33}) that \begin{equation}\label{eq:34}
\bigg| \frac{\mathcal{J}_{2,m}}{\mathcal{J}_{1,m}} \bigg| \leq
O\bigg( \bigg| \frac{a}{x - x_m} \bigg| \bigg) \leq O\bigg(
\frac{a}{d} \bigg) << 1.
\end{equation} From (\ref{eq:34}) by the arguments similar to the given
in \cite{R509} one obtains (\ref{eq:13}). \end{proof}
\begin{lmm}
Relation (\ref{eq:17}) holds. \end{lmm}
\begin{proof}
Let us justify relation (\ref{eq:17}). As $a \rightarrow 0$, one has
\begin{equation}\label{eq:35} \frac{\partial}{\partial
N_s}\frac{e^{-\sqrt{\lambda}|s - s'|}}{4\pi|s - s'|} =
\frac{\partial}{\partial N_s}\frac{1}{4\pi|s - s'|} +
\frac{\partial}{\partial N_s} \frac{e^{-\sqrt{\lambda}|s - s'|} -
1}{4\pi|s - s'|}. \end{equation} It is known (see \cite{R476}) that
\begin{equation}\label{eq:36}
\int_{\mathcal{S}_m}ds\int_{\mathcal{S}_m}\frac{\partial}{\partial
N_s}\frac{1}{4\pi|s - s'|}\sigma_m(s')ds' = -
\int_{\mathcal{S}_m}\sigma_m(s')ds' = -Q_m. \end{equation}
On the other hand, as $a \rightarrow 0$, one has
\begin{equation}\label{eq:37} \bigg|
\int_{\mathcal{S}_m}ds\int_{\mathcal{S}_m}\frac{e^{-\sqrt{\lambda}|s -
s'|} - 1}{4\pi|s - s'|}\sigma_m(s')ds' \bigg| \leq
|Q_m|\int_{\mathcal{S}_m}ds\frac{1-e^{-\sqrt{\lambda}|s - s'|}}{4\pi|s
- s'|} = o(Q_m). \end{equation}
The relations (\ref{eq:36}) and (\ref{eq:37}) justify (\ref{eq:17}).
\end{proof}
\begin{lmm}
Relation (\ref{eq:30}) holds. \end{lmm}
\begin{proof}
Denote
$$\displaystyle\frac{1}{t}\int^{t}_{0}u(t)dt := v(t),\quad \bar{u}(\sigma)
:= \displaystyle\int^{\infty}_{0}e^{-\sigma t}u(t)dt.$$ Then
$$\bar{v}(\lambda) =
\displaystyle\int^{\infty}_{\lambda}\frac{\bar{u}(\sigma)}{\sigma}d\sigma$$
by the properties of the Laplace transform. Assume that the limit
$v(\infty) := v_{\infty}$ exists: \begin{equation}\label{eq:38}
\lim_{t \rightarrow \infty}v(t) = v_{\infty}. \end{equation} Then,
$$v_{\infty} =
\displaystyle\lim_{\lambda \rightarrow
0}\lambda\displaystyle\int^{\infty}_{0}e^{-\lambda t}v(t)dt =
\displaystyle\lim_{\lambda \rightarrow 0}\lambda\bar{v}(\lambda).$$
Indeed,
$$\lambda\displaystyle\int^{\infty}_{0}e^{-\lambda t}dt =
1,$$
so
$$\displaystyle\lim_{\lambda \rightarrow
0}\lambda\displaystyle\int^{\infty}_{0}e^{-\lambda t}(v(t) -
v_{\infty})dt = 0,$$ and (\ref{eq:38}) is verified. One has
\begin{equation}\label{eq:39} \lim_{\lambda \rightarrow
0}\lambda\bar{v}(\lambda) = \lim_{\lambda \rightarrow
0}\int^{\infty}_{\lambda}\frac{\lambda}{\sigma}\bar{u}(\sigma)d\sigma
= \lim_{\lambda \rightarrow 0}\lambda\bar{u}(\lambda).
\end{equation} Let us check this:
\begin{equation}\label{eq:40} \lim_{\lambda \rightarrow
0}\int^{\infty}_{\lambda}\frac{\lambda}{\sigma}\bar{u}(\sigma)d\sigma =
\lim_{\lambda \rightarrow
0}\int^{\infty}_{\lambda}\frac{\lambda}{\sigma^2}\sigma\bar{u}(\sigma)d\sigma
= \lim_{\sigma \rightarrow 0}\sigma\bar{u}(\sigma), \end{equation} where
we have used the relation
$\displaystyle\int^{\infty}_{\lambda}\displaystyle\frac{\lambda}{\sigma^2}d\sigma
= 1$.$\\$ Alternatively, let $\sigma^{-1} = \gamma$. Then,
\begin{equation}\label{eq:41}
\int^{\infty}_{\lambda}\frac{\lambda}{\sigma^2}\sigma\bar{u}(\sigma)d\sigma
=
\frac{1}{1/\lambda}\int^{1/\lambda}_{0}\frac{1}{\gamma}\bar{u}(\frac{1}{\gamma})d\gamma
=
\frac{1}{\omega}\int^{\omega}_{0}\frac{1}{\gamma}\bar{u}(\frac{1}{\gamma})d\gamma.
\end{equation} If $\lambda \rightarrow 0$, then $\omega =
\lambda^{-1} \rightarrow \infty,$ and if $\psi :=
\gamma^{-1}\bar{u}(\gamma^{-1})$, then
\begin{equation}\label{eq:42} \lim_{\omega \rightarrow
\infty}\frac{1}{\omega}\int^{\omega}_{0}\psi d\gamma = \psi(\infty) =
\lim_{\gamma \rightarrow \infty}\gamma^{-1}\bar{u}(\gamma^{-1})
= \lim_{\sigma \rightarrow 0}\sigma\bar{u}(\sigma). \end{equation}
\end{proof}





\end{document}